\documentclass[a4paper,11pt,reqno]{amsart}
\usepackage{graphicx, psfrag, amsmath, amscd, amssymb}

\newtheorem{lem}{Lemma}[section]
\newtheorem{thm}[lem]{Theorem}
\newtheorem{pro}[lem]{Proposition}
\newtheorem{cor}[lem]{Corollary}
\newtheorem{exa}[lem]{Example}
\newtheorem{con}[lem]{Conjecture}
\newtheorem{defi}[lem]{Definition}

\newcommand{\ZZ}{{\mathbb{Z}}}

\newcommand{\U}{{\textsf{U}}}
\newcommand{\D}{{\textsf{D}}}

\title[Hankel Determinants for Dyck Paths]{On Hankel Determinants for Dyck Paths with Peaks Avoiding Multiple Classes of Heights}

 \author[H.-L. Chien]{Hsu-Lin Chien}
 \address{Department of Mathematics, National Taiwan Normal University, Taipei 116325, Taiwan, ROC}
 \email{ray010602@gmail.com}

 \author[{S.-P. Eu}]{Sen-Peng Eu}
 \address{Department of Mathematics, National Taiwan Normal University, Taipei 116325, and Chinese Air Force Academy, Kaohsiung 820009, Taiwan, ROC}
 \email{speu@math.ntnu.edu.tw}

 \author[T.-S. Fu]{Tung-Shan Fu}
 \address{Department of Applied Mathematics, National Pingtung University, Pingtung 900391, Taiwan, ROC}
 \email{tsfu@mail.nptu.edu.tw}

\begin{document}

\begin{abstract} For any integer $m\ge 2$ and a set $V\subset \{1,\dots,m\}$, let $(m,V)$ denote the union of congruence classes of the elements in $V$ modulo $m$. We study the Hankel determinants for the number of Dyck paths with peaks avoiding the heights in the set $(m,V)$. For any set $V$ of even elements of an even modulo $m$, we give an explicit description of the sequence of Hankel determinants in terms of subsequences of arithmetic progression of integers. There are numerous instances for varied $(m,V)$ with periodic sequences of Hankel determinants. We present a sufficient condition for the set $(m,V)$ such that the sequence of Hankel determinants is periodic, including even and odd modulus $m$.
\end{abstract}

\maketitle

\section{Introduction} 
\subsection{Hankel Matrices of generating functions}
Given a sequence $(f_n)_{n\ge 0}=(f_0,f_1,f_2,\dots)$, let  $F=F(x)=f_0+f_1 x+f_2 x^2+\cdots$ be the generating function of $(f_n)_{n\ge 0}$. Usually, $f_0=1$. For $n\ge 1$ and $k\ge 0$,  the  \emph{Hankel determinant}  of the series $F$ is defined by
\[
H^{(k)}_n(F):=\left|
              \begin{array}{cccc}
               f_k & f_{k+1} & \cdots & f_{k+n-1} \\
               f_{k+1} & f_{k+2} & \cdots & f_{k+n} \\
               \vdots  & \vdots  & \ddots & \vdots \\
               f_{k+n-1} & f_{k+n} & \cdots & f_{k+2n-2}      
              \end{array} 
             \right|
\]
By convention, we write $H_n(F)=H^{(0)}_n(F)$ for short. The \emph{sequence of Hankel determinants} of the series $F$ is defined to be
\[
H(F):=(H_1(F),H_2(F),H_3(F),\dots).
\]

There are a number of methods developed for evaluating Hankel determinants, such as continued fractions, orthogonal polynomials, $S$- and $J$-fractions (see \cite{Kr05,Viennot}), and $H$-fractions \cite{Han_fraction}. Recently, Han derived an explicit formula for the Hankel determinants of the Euler numbers by using $H$-fractions \cite{Han_Euler}.

Gessel--Viennot--Lindstr{\"o}m theorem is of combinatorial significance in linking the Hankel determinants of classical sequences to the enumeration of non-intersecting lattice paths \cite{GV}. For example, the Hankel matrices of Catalan numbers \cite{Aigner_ACCS,MW} (Motzkin numbers \cite{Viennot}, respectively) have determinant 1.
Non-intersecting Schr\"oder paths connect the Hankel determinant $2^{(n+1)n/2}$ of large Schr\"oder numbers to Aztec Diamond theorem \cite{BK,EF}. Other interesting sequences related to lattice paths include Catalan-like numbers \cite{Aigner} and weighted countings of partial Motzkin paths \cite{CK}.

Gessel and Xin  \cite{GX} developed a method of transformation, using bivariate polynomials, for the evaluation of Hankel determinants of a series. By Gessel--Xin's method, Sulanke and Xin  \cite{SX} used successive transformations to prove the periodicity of the sequence of Hankel determinants for the lattice paths with steps $\{(1,1),(3,0),(1,-1)\}$.

It seems that periodic sequences of Hankel determinants are rare. The notion of shifted periodic continued fractions was introduced in \cite{WXZ} and studied in \cite{WX} to evaluate shifted Hankel determinants for lattice paths.
In this paper, we study the Hankel determinants for the number of Dyck paths with peaks avoiding the set of heights consisting of multiple congruence classes of a modulo $m$. There are numerous instances for varied avoiding sets with periodic sequences of Hankel determinants. As a new technique, we establish fundamental reduction rules in recurrence form for evaluating the determinants, using Gessel--Xin's method of transformation.

\subsection{Dyck paths with restrictions on peaks}
A \emph{Dyck path} of \emph{size} $n$ is a lattice path from the origin to the point $(2n,0)$, using \emph{up-steps} $(1,1)$ and \emph{down-steps} $(1,-1)$, staying weakly above the $x$-axis. The number of Dyck paths of size $n$ is the $n$th Catalan number $\frac{1}{n+1}{\binom{2n}{n}}$. A \emph{peak} of a Dyck path is an up-step followed by a down-step. The \emph{height} of a peak is the $y$-coordinate of the intersection point of its steps. It is known that the Dyck paths with no peaks at odd (even, respectively) heights are counted by Riordan (shifted Motzkin, respectively) numbers \cite{ELY}.

Let $\ZZ^+$ denote the set of positive integers, and let $[n]:=\{1,2,\dots,n\}$ for any $n\in\ZZ^+$.
For any integer $m\ge 2$ and a set $V\subset [m]$, let $(m,V)$ denote the union of congruence classes of the elements in $V$ modulo $m$, namely,
\[
(m,V) := \{k\in \ZZ^+\, |\, k\equiv j\pmod{m} \mbox{ for some $j\in V$}\}.
\]
Moreover, let $V+t$ denote the set $\{ k\in [m]\, |\, k\equiv j+t\pmod{m} \mbox{ for some $j\in V$}\}$.

Let $d_n^{(m,V)}$ denote the number of Dyck paths $\pi$ of size $n$ such that the height of each peak of $\pi$ is not in the set $(m,V)$.
The generating function for the numbers $d_n^{(m,V)}$ is  
\[
D^{(m,V)}=D^{(m,V)}(x)=\sum_{n\ge 0} d_n^{(m,V)} x^n.
\]
The following relation between the series $D^{(m,V)}$ and $D^{(m,V-1)}$ is obtained from the `first return decomposition' of a Dyck path: $\pi=\U\mu\D\nu$, where $\mu$, $\nu$ are Dyck paths, and $\U$ ($\D$, respectively) is an up-step (down-step, respectively). Notice that $\mu$ is non-empty if $1\in V$.
\begin{equation} \label{eqn:CF_Dyck} 
D^{(m,V)}=\begin{cases}
           \dfrac{1}{1-xD^{(m,V-1)}} & \mbox{if $1\not\in V$} \\
           \dfrac{1}{1+x-xD^{(m,V-1)}} & \mbox{if $1\in V$.}
  \end{cases}.
\end{equation}
A closed form of the series $D^{(m,V)}$ can be obtained by solving the following equation, which is derived recursively from Eq.\,(\ref{eqn:CF_Dyck}); see \cite{ELY}. 
We make use of the $\chi$-notation that maps each statement $P$ onto $\{0, 1\}$, defined as $\chi(P)=1$ if $P$ is true, and is 0 otherwise. 
\begin{equation} \label{eqn:continued_fraction}
D^{(m,V)}=\dfrac{1}{1+\chi(1\in V)x-\dfrac{x}{\qquad\dfrac{\ddots}{1+\chi(m-1\in V)x-\dfrac{x}{1+\chi(m\in V)x-xD^{(m,V)}}}}}.
\end{equation}
For example, take $m=3$ and $V=\{1\}$. We obtain
\begin{align*}
D^{(3,\{1\})}(x) &= \frac{1-x^2-\sqrt{1-4x+2x^2+x^4}}{2x}\\
                 &= 1+x^2+2x^3+5x^4+13x^5+35x^6+97x^7+275x^8+794x^9+\cdots. 
\end{align*}

By Eq.\,(\ref{eqn:continued_fraction}), we remark that the series $D^{(m,V)}$ admits of an expansion in the form of Han's  \emph{super $\delta$-fraction} \cite{Han_fraction}, 
\[
D^{(m,V)}=\dfrac{v_0x^{k_0}}{1+u_1(x)x-\dfrac{v_1x^{k_0+k_1+\delta}}{1+u_2(x)x-\dfrac{v_2x^{k_1+k_2+\delta}}{1+u_3(x)-\ddots}}}
\]
with $\delta=1$, $v_j=1$, $k_j=0$, and $u_j(x)=\chi(j\in (m,V))$ for all $j$, but falls outside of the framework of super $\delta$-fraction that requests the degree of $u_j(x)$ less than or equal to $k_{j-1}+\delta-2$. There is a nice continued fraction for $D^{(m,V)}$ to be found.
Note that the classical methods of $S$- or $J$-fractions \cite{Kr99, Kr05, Viennot} do not apply to evaluating $H_n\big(D^{(m,V)}\big)$ because some of the Hankel determinants are zero.   

A \emph{periodic sequence} is written in contracted form using the notation with a star sign. Sometimes it is convenient to describe the structure of periodicity by using the form of an ultimately periodic sequence. For instance, the sequence $(1,(0,-1,1)^*)$ represents $(1$, $0,-1,1$, $0,-1,1$, $0,-1,1$,$\dots)=(1,0,-1)^*$. We observe that the sequence of Hankel determinants of the above series $D^{(3,\{1\})}$ is periodic with a period of 10 (Proposition \ref{pro:modulo_3}),
\[
H\big(D^{(3,\{1\})}\big) = (1, 1, 0, -1, -1, -1, -1, 0, 1, 1 )^*.
\]
Using the algorithm described in \cite{Han_fraction}, the periodic Hankel determinants $H(D^{(m,V)})$ in Example \ref{exa:example_1}(ii), Example \ref{exa:example_1.5}, Example \ref{exa:example_1.7}, and the periodicity in \cite{SX} can be proved automatically by computer. 

\subsection{Main results}
For the congruence of a modulo $\ell$,  let $\overline{k}$ denote the element congruent to $k$ mod $\ell$ with $1\le \overline{k}\le \ell$ for any positive integer $k$. 
Our first aim is to determine the sequence of Hankel determinants of the series $D^{(m,V)}$ for any set $V$ of even elements of an even modulo $m$.

\begin{thm} \label{thm:Main I} For any even integer $m\ge 2$ and a set $V=\{2s_1,\dots,2s_{\ell}\}$ with $1\le s_1<\cdots<s_{\ell}\le \frac{m}{2}$ and $\ell\ge 1$, let $t_j=s_{j+1}-s_j$ for $1\le j\le \ell-1$, and let $t_{\ell}=\frac{m}{2}-s_{\ell}+s_1$.
Then the sequence of Hankel determinants of the series $D^{(m,V)}$ can be partitioned into sections of the form
\begin{equation} \label{eqn:arithmetic_progression}
\begin{aligned}
H(D^{(m,V)}) &=(\underbrace{1,\dots,1}_{s_1},\underbrace{a_{1,1},\dots,a_{1,t_1}}_{t_1},\underbrace{a_{2,1},\dots,a_{2,t_2}}_{t_2},\dots,\underbrace{a_{\ell,1},\dots,a_{\ell,t_{\ell}}}_{t_{\ell}},  \\
             & \qquad \underbrace{a_{\ell+1,1},\dots,a_{\ell+1,t_1}}_{t_1},\underbrace{a_{\ell+2,1},\dots,a_{\ell+2,t_2}}_{t_2},\dots,\underbrace{a_{2\ell,1},\dots,a_{2\ell,t_{\ell}}}_{t_{\ell}},\dots )
\end{aligned}
\end{equation}
such that each section $(a_{j,1}, \dots,a_{j,t_{\overline{j}}})$, $\overline{j}\equiv j\pmod{\ell}$, is an arithmetic progression with the common difference $d_j$, where
$a_{1,1}=0$, $d_1=-1$, $a_{j,1}= d_{j-1}$, and $d_{j} = d_{j-1}-a_{j-1,t_{\overline{j-1}}}$ for $j\ge 2$.
\end{thm}


\begin{exa} \label{exa:example_1} {\rm (i) Take $m=10$ and $V=\{2,8\}$, i.e., $t_1=3$ and $t_2=2$. The sequence $H\big(D^{(10,\{2,8\})}\big)$ can be decomposed into sections of arithmetic progressions as follows.
\begin{align*}
H\big(D^{(10,\{2,8\})}\big) &=\big(1, (0, -1, -2), (-1, 0), (1, 2, 3), (1, -1), (-2, -3, -4), (-1, 2), \\ 
&\qquad (3, 4, 5), (1, -3), (-4, -5, -6), (-1, 4), (5, 6, 7), (1, -5), \dots\big).
\end{align*}

(ii) Take $m=10$ and $V=\{4,8\}$, i.e., $t_1=2$ and $t_2=3$. We observe that $H(D^{(10,\{4,8\})})$ is periodic with a period of 10,
\begin{align*}
H\big(D^{(10,\{4,8\})}\big) &=\big(1, 1, \big((0, -1), (-1, -1, -1), (0, 1), (1, 1, 1)\big)^*\big) \\
      &=(1, 1, 0, -1, -1, -1, -1, 0, 1,1)^*.
\end{align*}
}
\end{exa}

For the avoiding sets $(m,V)$ in which $V$ contains elements with mixed parities or the modulus $m$ is odd, there are series $D^{(m,V)}$ with periodic sequences $H\big(D^{(m,V)}\big)$, but the others, with non-periodic sequences $H\big(D^{(m,V)}\big)$, are not expected to have nice results. For example,
\begin{align*}
H\big(D^{(5,\{2\})} \big) &=(1, 0, -1, -2, -2, -3, -4, -5, -1, 7, 23, 31, 51, 116, 149, 118, -426, \dots), \\
H\big(D^{(7,\{2\})} \big) &=(1, 0, -1, -2, -3, -3, -4, -8, -9, -10, -4, 26, 53, 104, 212, 323, 671, \dots), \\
H\big(D^{(5,\{1,2\})}\big) &=(1,0,-1,-1,-1,-1,0,1,1,1)^*, \\
H\big(D^{(6,\{1,2\})}\big) 
&=(1,0,-1,-2,0,2,5,8,11,3,3,-17,-260,-452, -839, -1752, 5288,\dots).
\end{align*}

Our second aim is to present a sufficient condition for the set $(m,V)$ such that the sequence $H\big(D^{(m,V)}\big)$ is periodic, including even and odd modulus $m$.

\begin{defi} {\rm
For $b\ge 0$, a sequence $(t_1,\dots,t_b)$ of positive integers is called \emph{admissible} if the sequence $\big(h(t_1),h(t_1,t_2),\dots,h(t_1,\dots,t_b)\big)$ of rational numbers, defined inductively by $h(t_1)=2-t_1$ and
\[ 
h(t_1,\dots,t_j) := 2-t_j-\frac{1}{h(t_1,\dots,t_{j-1})}=2-t_j-\dfrac{1}{2-t_{j-1}-\dfrac{1}{\qquad\dfrac{\ddots}{2-t_2-\dfrac{1}{2-t_1}}}}
\]
for $2\le j\le b$, are all non-zero. We assume that the empty sequence is admissible. Moreover, if there exists a positive integer $t_{b+1}$ such that $h(t_1,\dots,t_{b+1}):= 2-t_{b+1}-{h(t_1,\dots,t_b)}^{-1}$ equals zero then the sequence $(t_1,\dots,t_{b+1})$ is called \emph{primitive}.  
The sequence $\big(h(t_1),\dots,h(t_1,\dots,t_{b+1})\big)$ is called the \emph{dual sequence} associated with $(t_1,\dots,t_{b+1})$.
}
\end{defi}

For example, $(t_1)=(2)$ is the primitive sequence of length 1, $(t_1,t_2)=(1,1)$, $(3,3)$ are the primitive sequences of length 2, and $(t_1,t_2,t_3)=(1,2,3)$, $(3,2,1)$, $(3,4,3)$, and $(4,3,4)$ are the primitive sequences of length 3.

\begin{thm} \label{thm:Main II} For positive integers $t_1,\dots,t_{\ell-1}$ and $s$  ($\ell\ge 2$), let $V=\{2s,2(s+t_1),\dots,2(s+t_1+\cdots+t_{\ell-1})\}$ and $m\ge\max(V)$. If the sequence $(t_1,\dots,t_{\ell-1})$ is primitive then the sequence of Hankel determinants of the series $D^{(m,V)}$ is periodic with a period of either $p$ or $2p$, where
\[
p=\begin{cases}
\frac{m}{2} &\mbox{for $m$ even} \\
m &\mbox{for $m$ odd.}
\end{cases}
\]
Moreover, the periodicity is $p$ ($2p$, respectively) if the partial product $\prod_{j=1}^{\ell-2} h(t_1,\dots,t_j)$ of the dual sequence of $(t_1,\dots,t_{\ell-1})$ is  $-1$ ($1$, respectively).
\end{thm}

\begin{exa} \label{exa:example_1.5} {\rm  (i) As shown in Example \ref{exa:example_1}(ii), the sequence $H\big(D^{(10,\{4,8\})}\big)$ is periodic with a period of 10. Notice that  $V=\{4,8\}$ is one of the sets derived from the primitive sequence of length 1 (i.e., $t_1=2$).

(ii) Take $m=16$ and $V=\{4,10,16\}$, a set derived from the primitive sequence $(3,3)$. We observe that the periodicity of the sequence $H\big(D^{(16,\{4,10,16\})}\big)$ is 8,
\begin{align*}
H\big(D^{(16,\{4,10,16\})}\big) 
&=\big(1,1,((0,-1,-2),(-1,0,1),(1,1))^*   \big) \\
&= (1, 1, 0, -1, -2, -1, 0, 1)^*.
\end{align*}

(iii) Take $m=17$ and $V=\{4,10,16\}$.  The periodicity of the sequence $H\big(D^{(17,\{4,10,16\})}\big)$ is 17,
\begin{align*}
H\big(D^{(17,\{4,10,16\})}\big) 
&= \big(1, 1, ((0, -1, -2), (-1, 0, 1), (1, 1, 1, 1, 1, 1, 1, 1, 1,1,1))^*  \big) \\
&= (1, 1, 0, -1, -2, -1, 0, 1, 1, 1, 1, 1, 1, 1, 1, 1, 1)^*.
\end{align*}

(iv) Take $m=10$ and $V=\{2,4,6\}$, a set derived from the primitive sequence $(1,1)$. We observe that the periodicity of the sequence $H\big(D^{(10,\{2,4,6\})}\big)$ is 10,
\begin{align*}
H\big(D^{(10,\{2,4,6\})}\big) 
&= \big(1, ((0), (-1), (-1, -1, -1), (0), (1), (1, 1, 1))^* \big) \\
&= (1, 0, -1, -1, -1, -1, 0, 1, 1, 1)^*.
\end{align*}

(v) Take $m=11$ and $V=\{2,4,6\}$. The periodicity of the sequence $H\big(D^{(11,\{2,4,6\})}\big)$ is 22,
\begin{align*}
H\big(D^{(11,\{2,4,6\})}\big)  
&= \big(1, ((0), (-1), (-1, -1, -1, -1, -1, -1, -1, -1, -1), (0), (1),\\
&\qquad  (1, 1, 1, 1, 1, 1, 1, 1, 1))^* \big) \\
&= (1, 0, -1, -1, -1, -1, -1, -1, -1, -1, -1, -1, 0, 1, 1, 1, 1, 1, 1, 1, 1, 1)^*.
\end{align*}
}
\end{exa}

A set $V\subset Z^{+}$ of elements is called \emph{feasible} if the sequence $H\big(D^{(m,V)}\big)$ is periodic for some integer $m$. As shown in Theorem \ref{thm:Main II}, a feasible set derived from a primitive sequence $\alpha=(t_1,\dots,t_{\ell-1})$ is of the form $V_{\alpha}=\{2s,2(s+t_1),\dots,2(s+t_1+\cdots+t_{\ell-1})\}$ for some $s\ge 1$.
We construct a variety of feasible sets by a synthesis of the ones derived from primitive sequences.


\begin{thm} \label{thm:variety} Let $\alpha_1,\dots,\alpha_b$ be primitive sequences, and let $V_j$ be the feasible set derived from $\alpha_i$ with $\min(V_j)=2$ for $j=1,\dots, b$. If $V=(V_1+k_1)\cup\cdots\cup(V_b+k_b)$ for some $k_1\ge -1$, and $k_{j+1}\ge\max(V_j+k_j)-1$, $1\le j\le b-1$, then the sequence of Hankel determinants of the series $D^{(m,V)}$ is periodic for every $m\ge \max(V)$.
\end{thm}

\begin{exa} \label{exa:example_1.7} {\rm (i) Take $V_{\alpha}=\{2,6\}$, a feasible set derived from the primitive sequence of length 1. Let $V=V_{\alpha}\cup(V_{\alpha}+5)=\{2,6,7,11\}$. For $m=14$ ($m=11$, respectively), we have the following periodic sequences,  
\begin{align*}
H\big(D^{(14,\{2,6,7,11\})}\big) 
&= (1, 0, -1, -1, -1, -1, -1, -1, 0, 1, 1, 1, 1, 1)^*, \\
H\big(D^{(11,\{2,6,7,11\})}\big) 
&= (1, 0, -1, -1, -1, -1, -1, -1, -1, 0, 1)^*.
\end{align*}

(ii) If $V=(V_{\alpha}-1)\cup(V_{\alpha}+6)=\{1,5,8,12\}$ and $m=14$,  we have
\[
H\big(D^{(14,\{1,5,8,12\})}\big) 
= (1, 1, 1, 1, 0, -1, -1, -1, -1, -1, -1, 0, 1, 1)^*.
\]
Moreover, if $V=(V_{\alpha}-1)\cup(V_{\alpha}+4)=\{1,5,6,10\}$ and $m=11$,  we have
\[
H\big(D^{(11,\{1,5,6,10\})}\big) 
= (1, 1, 1, 0, -1, -1, 0, 1, 1, 1, 1)^*.
\]
}
\end{exa}

The rest of this paper is organized as follows. We review Gessel--Xin's method of transformation and establish fundamental reduction rules for the Hankel determinants $H_n\big(D^{(m,V)}\big)$ in Section 2. The proofs of our main results (Theorem \ref{thm:Main I} and Theorem \ref{thm:Main II}) are given in Section 3 and Section 4, respectively. A synthesis of feasible sets for periodic sequence of Hankel determinants (Theorem \ref{thm:variety}) is proved in Section 5. Some sporadic cases  not covered by the sufficient condition for periodicity are given in Section 6. For readers' reference, we list the periodicity of the sequence $H\big(D^{(m,V)}\big)$ for all $(m,V)$ with $m\le 5$ in Table \ref{tab:list} at the end of this paper.

\section{Reduction rules for the Hankel determinants $H_n\big(D^{(m,V)}\big)$}
In this section, we prove the following relations for the Hankel determinants of the series $D^{(m,V)}$ and $-1+D^{(m,V)}$, which pave the way to the proof of Theorem \ref{thm:Main I}.

\begin{pro} \label{pro:reduction_rules} For any integer $m\ge 2$ and a set $V\subset [m]$,  the Hankel determinants of the series $D^{(m,V)}$ and $-1+D^{(m,V)}$ satisfy the following relations. 
\begin{enumerate}
\item We have
\[
H_n\big(D^{(m,V)}\big)=\begin{cases}
          H_{n-1}\big(D^{(m,V-2)}\big) & \mbox{if $2\not\in V$} \\
          H_{n-1}\big(-1+D^{(m,V-2)}\big) & \mbox{if $2\in V$.}
  \end{cases}
\]
\item If $1\not\in V$ then we have
\[
H_n\big(-1+D^{(m,V)}\big)= 
-H_{n-1}\big(D^{(m,V-2)}\big)+H_{n-1}\big(-1+D^{(m,V-2)}\big).
\]
\end{enumerate}
\end{pro}

\medskip
The proof of Proposition \ref{pro:reduction_rules}(i) (Proposition \ref{pro:reduction_rules}(ii), respectively) is given in Lemmas \ref{lem:H(0)} and \ref{lem:H(1)} (Lemmas \ref{lem:-H(2)} and  \ref{lem:H(2)}, respectively). 

We remark that the pivotal points in the evaluation of $H_n\big(D^{(m,V)}\big)$, using Proposition \ref{pro:reduction_rules}, are the reductions when  the element 2 appears in the avoiding set.
There is an obstacle in the process when $1\in V-2j$ for some $j$ since no deterministic rule for $H_{n-j}\big(-1+D^{(m,V-2j)}\big)$ is available as in Proposition \ref{pro:reduction_rules}(ii). This is in connection with the chaos of non-periodic sequences $H\big(D^{(m,V)}\big)$ when $V$ contains elements with mixed parities or the modulus $m$ is odd.

\subsection{Manipulations for Hankel determinants}

\begin{lem} \label{lem:H(2)} For a series $F=F(x)=f_0+f_1x+f_2x^2+\cdots$, we have
\[
H^{(2)}_{n-1}(F)=-H_n(-1+F)+H_n(F).
\]
\end{lem}

\begin{proof}
\begin{align*}
H_n(F) &=\left|
              \begin{array}{cccc}
               f_0 & f_1 & \cdots & f_{n-1} \\
               f_1 & f_2 & \cdots & f_n \\
               \vdots  & \vdots  & \ddots & \vdots \\
               f_{n-1} & f_n & \cdots & f_{2n-2}      
              \end{array} 
             \right| \\
       &=\left|
              \begin{array}{cccc}
               f_0-1 & f_1 & \cdots & f_{n-1} \\
               f_1 & f_2 & \cdots & f_n \\
               \vdots  & \vdots  & \ddots & \vdots \\
               f_{n-1} & f_n & \cdots & f_{2n-2}      
              \end{array} 
             \right|+
             \left|
              \begin{array}{cccc}
               1 & f_1 & \cdots & f_{n-1} \\
               0 & f_2 & \cdots & f_n \\
               \vdots  & \vdots  & \ddots & \vdots \\
               0 & f_n & \cdots & f_{2n-2}      
              \end{array} 
             \right|. \\
       &= H_n(-1+F)+H_{n-1}^{(2)}(F).      
\end{align*}
The assertion follows.
\end{proof}

We shall evaluate $H_n(F)$ by successive reductions to $n=1$ and use the following facts.
\begin{equation} \label{eqn:trivial}
H_1(F)=f_0=1, \qquad\mbox{and}\qquad H_1(-1+F)=f_0-1=0.
\end{equation}

We review Gessel--Xin's method of transformation for  Hankel determinants of a series $F$ \cite{GX}. For a bivariate polynomial $G(x,y)=g_{i,j}x^i y^j$, let $[G(x,y)]_n$ denote the determinant of the $n\times n$ matrix $(g_{i,j})_{0\le i,j\le n-1}$. Hence the Hankel determinants $H_n(F)$ and $H^{(1)}_n(F)$ of the series $F$  can be expressed as
\begin{equation} \label{eqn:bivariate}
H_n(F)=\left[\frac{xF(x)-yF(y)}{x-y} \right]_n, \qquad H_n^{(1)}(F)=\left[\frac{F(x)-F(y)}{x-y}  \right]_n.
\end{equation}
The following rule for transformations of the determinant $[G(x,y)]_n$ holds.


\begin{lem} \label{lem:production_rule} {\rm (Product Rule)}  If $u(x)$ is a formal power series with $u(0)=1$ then
\[
[u(x)G(x,y)]_n=[u(y)G(x,y)]_n=[G(x,y)]_n.
\]
\end{lem}


\subsection{Proof of the reduction rules} We deduce the results in Lemmas \ref{lem:H(0)}, \ref{lem:H(1)} and \ref{lem:-H(2)} by applying Product Rule to the transformations of the series $D^{(m,V)}$.

\begin{lem} \label{lem:H(0)} For any integer $m\ge 2$ and a set $V\subset [m]$,  the following relation holds.
\[
H_n\big(D^{(m,V)}\big)=H_{n-1}^{(1)}\big(D^{(m,V-1)}\big).
\]
\end{lem}

\begin{proof} 
By Eq.\,(\ref{eqn:CF_Dyck}), if $1\not\in V$ then $D^{(m,V)}=\big(1-xD^{(m,V-1)}\big)^{-1}$. By the transformation in Eq.\,(\ref{eqn:bivariate}),  we have
\begin{equation*} \label{eqn:1_not_in_V}
  \begin{aligned}
H_n\big(D^{(m,V)}\big) &= \left[\dfrac{xD^{(m,V)}(x)-yD^{(m,V)}(y)}{x-y}\right]_n \\
 &= \left[\dfrac{\dfrac{x}{1-x D^{(m,V-1)}(x)}-\dfrac{y}{1-y D^{(m,V-1)}(y)}}{x-y}\right]_n
    \end{aligned}
\end{equation*}
By Product Rule, multiplying by the series $\big(1-xD^{(m,V-1)}(x)\big)\big(1-yD^{(m,V-1)}(y)\big)$ yields
\begin{equation*} \label{eqn:1_not_in_V_part2}
  \begin{aligned}
H_n\big(D^{(m,V)}\big) &= \left[\dfrac{x\big(1-yD^{(m,V-1)}(y)\big)-y\big(1-xD^{(m,V-1)}(x)\big)}{x-y}\right]_n \\
 &= \left[1+xy\dfrac{D^{(m,V-1)}(x)-D^{(m,V-1)}(y)}{x-y}\right]_n \\
 &= H^{(1)}_{n-1}\big(D^{(m,V-1)}\big).
   \end{aligned}
\end{equation*}

For the case $1\in V$, we have $D^{(m,V)}=\big(1+x-xD^{(m,V-1)}\big)^{-1}$. Making use of this relation in the evaluation of $H_n\big(D^{(m,V)}\big)$ as above, the assertion $H_n\big(D^{(m,V)}\big)=H^{(1)}_{n-1}\big(D^{(m,V-1)}\big)$ is also obtained.
\end{proof}

\begin{lem} \label{lem:H(1)} For any integer $m\ge 2$ and a set $V\subset [m]$,  the following relation holds.
\[
H^{(1)}_n(D^{(m,V)}\big)=\begin{cases}
          H_n\big(D^{(m,V-1)}\big) & \mbox{if $1\not\in V$} \\
          H_n\big(-1+D^{(m,V-1)}\big) & \mbox{if $1\in V$.}
  \end{cases}
\]
\end{lem}

\begin{proof}  (i) If $1\not\in V$ then $D^{(m,V)}=\big(1-xD^{(m,V-1)}\big)^{-1}$. By Eq.\,(\ref{eqn:bivariate}) and Lemma \ref{lem:production_rule}, we have
\begin{equation} \label{eqn:H(1)_1_not_in_V}
  \begin{aligned}
H^{(1)}_n\big(D^{(m,V)}\big) &= \left[\dfrac{D^{(m,V)}(x)-D^{(m,V)}(y)}{x-y}\right]_n \\
 &= \left[\dfrac{\dfrac{1}{1-x D^{(m,V-1)}(x)}-\dfrac{1}{1-y D^{(m,V-1)}(y))}}{x-y}\right]_n  \\
 &= \left[\dfrac{\big(1-yD^{(m,V-1)}(y)\big)-\big(1-xD^{(m,V-1)}(x)\big)}{x-y}\right]_n \\
 &= \left[\dfrac{xD^{(m,V-1)}(x)-yD^{(m,V-1)}(y)}{x-y}\right]_n \\ \\
 &= H_n\big(D^{(m,V-1)}\big).
   \end{aligned}
\end{equation}

(ii) If $1\in V$ then $D^{(m,V)}=\big(1+x-xD^{(m,V-1)}\big)^{-1}$. We have
\begin{align*}
H^{(1)}_n\big(D^{(m,V)}\big) &= \left[\dfrac{D^{(m,V)}(x)-D^{(m,V)}(y)}{x-y}\right]_n \\
 &= \left[\dfrac{\dfrac{1}{1+x-x D^{(m,V-1)}(x)}-\dfrac{1}{1+y-y D^{(m,V-1)}(y)}}{x-y}\right]_n  \\
 &= \left[\dfrac{\big(1+y-yD^{(m,V-1)}(y)\big)-\big(1+x-xD^{(m,V-1)}(x)\big)}{x-y}\right]_n \\
 &= \left[\dfrac{x\big(-1+D^{(m,V-1)}(x)\big)-y\big(-1+D^{(m,V-1)}(y)\big)}{x-y}\right]_n \\ 
 &= H_n\big(-1+D^{(m,V-1)}\big).
\end{align*}
\end{proof}

By Lemmas \ref{lem:H(0)} and \ref{lem:H(1)}, we prove Proposition \ref{pro:reduction_rules}(i).

\begin{lem} \label{lem:-H(2)} For any integer $m\ge 2$ and a set $V\subset [m]$, if $1\not\in V$ then the following relation holds.
\[
H_n\big(-1+D^{(m,V)}\big)=-H^{(2)}_{n-2}\big(D^{(m,V-2)}\big).
\]
\end{lem}

\begin{proof} Notice that $D^{(m,V)}=\big(1-xD^{(m,V-1)}\big)^{-1}$ since $1\not\in V$.  There are two cases.

Case 1. If $2\in V$ then $1\in V-1$ and $D^{(m,V-1)}=\big(1+x-xD^{(m,V-2)}\big)^{-1}$. We have
\begin{equation} \label{eqn:2_in_V}
  \begin{aligned}
-1+D^{(m,V)} &=-1+\frac{1}{1-xD^{(m,V-1)}} \\
                &=\frac{xD^{(m,V-1)}}{1-xD^{(m,V-1)}} \\
                &=\frac{x\big(1+x-xD^{(m,V-2)}\big)^{-1}}{1-x\big(1+x-xD^{(m,V-2)}\big)^{-1}} \\
                &=\frac{x}{1-xD^{(m,V-2)}}.
   \end{aligned}
\end{equation}
By Eq.\,(\ref{eqn:bivariate}) and the above relation, the Hankel determinant of the series $-1+D^{(m,V)}$ is
\begin{equation} \label{eqn:2x2_determinant}
  \begin{aligned}
H_n\big(-1+D^{(m,V)}\big) &= \left[\dfrac{x\big(-1+D^{(m,V)}(x)\big)-y\big(-1+D^{(m,V)}(y)\big)}{x-y}\right]_n \\
 &= \left[\dfrac{\dfrac{x^2}{1-xD^{(m,V-2)}(x)}-\dfrac{y^2}{1-y D^{(m,V-2)}(y)}}{x-y}\right]_n  \\
 &= \left[\dfrac{x^2\big(1-yD^{(m,V-2)}(y)\big)-y^2\big(1-xD^{(m,V-2)}(x)\big)}{x-y}\right]_n \\
 &= \left[\frac{1}{x-y}\Bigg( (x^2-y^2)-xy(x-y)  \right.\\
 &\qquad\qquad \left.  +x^2y^2\left(\frac{-1+D^{(m,V-2)}(x)}{x}-\frac{-1+D^{(m,V-2)}(y)}{y}\right)\Bigg)\right]_n \\ \\
 &= \left| \begin{array}{cc}
                0 & 1 \\
                1 & -1
           \end{array}
    \right| \cdot
 H^{(2)}_{n-2}\big(D^{(m,V-2)}\big).
   \end{aligned}
\end{equation}

Case 2. If $2\not\in V$ then $1\not\in V-1$ and  $D^{(m,V-1)}=\big(1-xD^{(m,V-2)}\big)^{-1}$. 
By the same argument as in Eq.\,(\ref{eqn:2_in_V}), we have
\begin{equation} \label{eqn:2_not_in_V}
-1+D^{(m,V)}=\frac{x}{1-x-xD^{(m,V-2)}}.
\end{equation}
Making use of Eq.\,(\ref{eqn:2_not_in_V}) in the evaluation of $H_n\big(-1+D^{(m,V)}\big)$ as in Eq.\,(\ref{eqn:2x2_determinant}), we obtain
\[
H_n\big(-1+D^{(m,V)}\big)=\left| \begin{array}{cc}
                0 & 1 \\
                1 & -2
           \end{array}
    \right| \cdot
 H^{(2)}_{n-2}\big(D^{(m,V-2)}\big).
\]
The assertion follows from both of the cases.
\end{proof}

By Lemmas \ref{lem:-H(2)} and \ref{lem:H(2)}, we prove Proposition \ref{pro:reduction_rules}(ii). 
%
%

\section{Proof of Theorem \ref{thm:Main I}} In this section, we prove Theorem \ref{thm:Main I}. First, we present a proof for the case when the set $V$ is a singleton. In particular, this proves Viennot's result \cite{Viennot} on the periodicity of the sequence of Hankel determinants of shifted Motzkin numbers (see Corollary \ref{cor:Viennot_result}).

\begin{thm} \label{thm:singleton} If $m=2p\ge 2$ and $V=\{2s\}$ for some $s$, $1\le s\le p$, then the sequence of Hankel determinants of the series $D^{(m,V)}$ can be partitioned into sections of the form
\[
H\big(D^{(2p,\{2s\})}\big)=(\underbrace{1,\dots,1}_{s},\underbrace{a_{1,1},\dots,a_{1,p}}_{p},\underbrace{a_{2,1},\dots,a_{2,p}}_{p},\underbrace{a_{3,1},\dots,a_{3,p}}_{p},\dots)
\]
such that each section $(a_{j,1}, \dots,a_{j,p})$ is an arithmetic progression with the common difference $d_j$, where $a_{1,1}=0$, $d_1= -1$, $a_{j,1}= d_{j-1}$,  and $d_j = d_{j-1}-a_{j-1,p}$ for $j\ge 2$.
\end{thm}

\begin{proof} We shall apply the reduction rules in Proposition \ref{pro:reduction_rules} to the evaluation of $H_n\big(D^{(m,V)}\big)$ successively. Notice that $1\not\in V-2j$, modulo $m$, for all $j\ge 0$, and that the pivotal points are reductions at $V-mj$, containing the element 2.

(i) Let $s=1$. For $1\le k\le p$, by Proposition \ref{pro:reduction_rules} and Eq.\,(\ref{eqn:trivial}), we have the initial elements
\begin{equation} \label{eqn:initial_section}
\begin{aligned}
a_{1,k}&=H_{k+1}\big(D^{(m,\{2\})} \big) \\
      &=H_k\big(-1+D^{(m,\{m\})} \big) \\  
      &=-H_{k-1}\big(D^{(m,\{m-2\})}\big)+H_{k-1}\big(-1+D^{(m,\{m-2\})}\big) \\
      &=\dots \\
      &=-(k-1)H_1\big(D^{(m,\{m-2k+2\})}\big)+H_1\big(-1+D^{(m,\{m-2k+2\})}\big) \\
      &=-k+1.    
\end{aligned}
\end{equation}
Hence $a_{1,1}=0$ and $d_1=-1$. For $i\ge 1$, let $d_{i+1}=d_i-a_{i,p}$. We shall prove that $a_{i+1, k}=d_i+(k-1)d_{i+1}$ for $1\le k\le p$. It suffices to prove the following relation
\begin{equation}  \label{eqn:Induction}
H_{ip+t}\big(D^{(m,\{2\})}\big)=a_{i,p}H_t\big(D^{(m,\{2\})}\big)-d_iH_t\big(-1+D^{(m,\{2\})}\big), \quad t\ge 2.
\end{equation}

For $i=1$, by the argument in Eq.\,(\ref{eqn:initial_section}), we observe that
\begin{align*}
H_{p+t}\big(D^{(m,\{2\})} \big)
      &=H_{p+t-1}\big(-1+D^{(m,\{m\})} \big) \\
      &=\cdots \\
      &=-(p-1)H_t\big(D^{(m,\{2\})}\big)+H_t\big(-1+D^{(m,\{2\})}\big) \\
      &=a_{1,p}H_t\big(D^{(m,\{2\})}\big)-d_1H_t\big(-1+D^{(m,\{2\})}\big).
\end{align*}
Suppose the relation Eq.\,(\ref{eqn:Induction}) holds for $i\le j$. By Proposition \ref{pro:reduction_rules} and Eq.\,(\ref{eqn:trivial}), we have 
\begin{equation} \label{eqn:a_j+1_1}
\begin{aligned}
a_{j+1,1}&=H_{jp+2}\big(D^{(m,\{2\})} \big) \\
      &=a_{j,p}H_{2}\big(D^{(m,\{2\})} \big)-d_jH_2\big(D^{(m,\{2\})} \big)\\
      &=d_jH_1\big(D^{(m,\{m\})}\big)+(a_{j,p}-d_j)H_1\big(-1+D^{(m,\{m\})}\big) \\
      &=d_j.
\end{aligned}    
\end{equation}  
By induction hypothesis, we have 
\begin{equation} \label{eqn:induction_hypothesis}
\begin{aligned}       
H_{(j+1)p+t}\big(D^{(m,\{2\})} \big)
      &=a_{j,p}H_{p+t}\big(D^{(m,\{2\})}\big)-d_jH_{p+t}\big(-1+D^{(m,\{2\})}\big) \\ 
      &=d_jH_{p+t-1}\big(D^{(m,\{m\})}\big)+(a_{j,p}-d_j)H_{p+t-1}\big(-1+D^{(m,\{m\})}\big) \\
      &=a_{j+1,1}H_{p+t-1}\big(D^{(m,\{m\})}\big)-d_{j+1}H_{p+t-1}\big(-1+D^{(m,\{m\})}\big) \\
      &=\cdots \\
      &=(a_{j+1,1}-(p-1)d_{j+1})H_t\big(D^{(m,\{2\})}\big)-d_{j+1}H_t\big(-1+D^{(m,\{2\})}\big) \\
      &=a_{j+1,p}H_t\big(D^{(m,\{2\})}\big)-d_{j+1}H_t\big(-1+D^{(m,\{2\})}\big).
\end{aligned}
\end{equation}
The result $a_{i+1, k}=d_i+(k-1)d_{i+1}$ can be derived from the relation Eq.\,(\ref{eqn:Induction}), by the same argument as in Eqs.\,(\ref{eqn:a_j+1_1}) and (\ref{eqn:induction_hypothesis}).

(ii) If $s>1$ then for $n\ge 1$ we have 
\[
H_{n+s}\big(D^{(m,\{2s\})} \big)=H_{n+s-1}\big(D^{(m,\{2s-2\})} \big)=\cdots=H_{n+1}\big(D^{(m,\{2\})} \big).
\]
For $1\le k\le s$, notice that $H_k\big(D^{(m,\{2s\})})=1$ is the Hankel determinant of Catalan numbers.
\end{proof}

\begin{exa} {\rm Take $m=10$ and $V=\{2\}$. The sequence $H\big(D^{(10,\{2\})}\big)$ can be decomposed into sections of arithmetic progressions as follows.
\begin{align*}
H\big(D^{(10,\{2\})}\big) &= \big(1, (0,-1,-2,-3,-4), (-1,2,5,8,11), (3,-5,-13,-21,-29),\\
&\qquad (-8,13,34,55,76), (21,-34,-89,-144,-199), (-55,89,233,377,521), \dots\big). \\ 
\end{align*}
}
\end{exa}

By Theorem \ref{thm:singleton}, we have the following  immediate result.

\begin{cor} \label{cor:Viennot_result}
The sequence of Hankel determinants for the Dyck paths with no peaks at even heights is periodic with a period of 6, namely,
\[
H\big(D^{(2,\{2\})}\big)=(1,0,-1,-1,0,1)^*.
\]
\end{cor}

Now, we prove Theorem \ref{thm:Main I} by extending the proof of Theorem \ref{thm:singleton}.

\medskip
\noindent
\emph{Proof of Theorem \ref{thm:Main I}.}
(i)  Let $s=1$, i.e., $V=\{2,2(1+t_1),\dots,2(1+t_1+\cdots+t_{\ell-1})\}$. Consider the decomposition of $H\big(D^{(m,V)}\big)$ in Eq.\,(\ref{eqn:arithmetic_progression}). For the initial elements, we have $a_{1,k}=-k+1$ for $1\le k\le t_1$ by the same argument as in Eq.\,(\ref{eqn:initial_section}). Moreover, for $1\le k\le t_2$, we have
\begin{equation} \label{eqn:initial_item_section}
\begin{aligned}
a_{2,k}&=H_{t_1+k+1}\big(D^{(m,V)} \big) \\
      &=H_{t_1+k}\big(-1+D^{(m,V-2)} \big) \\  
      &=\dots \\
      &=-(t_1-1)H_{k+1}\big(D^{(m,V-2t_1)}\big)+H_{k+1}\big(-1+D^{(m,V-2t_1)}\big) \\
      &=a_{1,t_1}H_{k+1}\big(D^{(m,V-2t_1)}\big)-d_1H_{k+1}\big(-1+D^{(m,V-2t_1)}\big).
\end{aligned}
\end{equation}
Notice that $2\in V-2(t_1+\cdots+t_i)$ for $1\le i\le \ell-1$ and that $V-2(t_1+\cdots+t_{\ell})=V-m=V$. For $1\le b\le \ell$, let $d_{b+1}=d_b-a_{b,t_b}$. By the same argument as in Eq.\,(\ref{eqn:initial_item_section}), we observe that
\begin{equation} \label{eqn:first_run}
\begin{aligned}
a_{b+1,k} &=H_{t_1+\cdots+t_b+k+1}\big(D^{(m,V)} \big) \\
       &=a_{1,t_1}H_{t_2+\cdots+t_b+k+1}\big(D^{(m,V-2t_1)} \big)-d_1H_{t_2+\cdots+t_b+k+1}\big(-1+D^{(m,V-2t_1)}\big)\\
       &=\cdots \\
       &=a_{b,t_b}H_{k+1}\big(D^{(m,V-2(t_1+\cdots+t_b))} \big)-d_bH_{k+1}\big(-1+D^{(m,V-2(t_1+\cdots+t_b))}\big) \\
       &=d_b+(k-1)d_{b+1}   
\end{aligned}
\end{equation}
for $1\le k\le t_{b+1}$. Let $p=\frac{m}{2}$. Notice that $2\in V-mj=V$ for all $j$.  For $i\ge 1$ and $1\le b\le \ell$, let $d_{i\ell+b+1}=d_{i\ell+b}-a_{i\ell+b,t_b}$. By the same argument as in Eqs.\,(\ref{eqn:induction_hypothesis}) and (\ref{eqn:first_run}), we have
\begin{equation} \label{eqn:general_run}
\begin{aligned}
a_{i\ell+b+1,k} &=H_{ip+t_1+\cdots+t_b+k+1}\big(D^{(m,V)} \big) \\
       &=a_{\ell,t_{\ell}}H_{(i-1)p+t_1+\cdots+t_b+k+1}\big(D^{(m,V)} \big)-d_{\ell}H_{(i-1)p+t_1+\cdots+t_b+k+1}\big(-1+D^{(m,V)}\big)\\
       &=\cdots \\
       &=a_{i\ell,t_{\ell}}H_{t_1+\cdots+t_b+k+1}\big(D^{(m,V)} \big)-d_{i\ell}H_{t_1+\cdots+t_b+k+1}\big(-1+D^{(m,V)}\big) \\
       &=d_{i\ell+b}+(k-1)d_{i\ell+b+1}.   
\end{aligned}
\end{equation}

(ii) For $V=\{2s,2(s+t_1),\dots,2(s+t_1+\cdots+t_{\ell-1})\}$ and $s>1$, notice that $2\in V-2(s-1)$. Then for $n\ge 1$, we have 
\[
H_{n+s}\big(D^{(m,V)} \big)=H_{n+s-1}\big(D^{(m,V-2)} \big)=\cdots=H_{n+1}\big(D^{(m,V-2(s-1))} \big).
\]
For $1\le k\le s$, notice that $H_k\big(D^{(m,V)})=1$ is the Hankel determinant of Catalan numbers.
\qed

\begin{exa} {\rm Take $m=24$ and $V=\{2,8,12,18\}$. The sequence $H\big(D^{(24,\{2,8,12,18\})}\big)$ can be decomposed into sections of arithmetic progressions as follows.
\begin{align*}
H\big(D^{(24,\{2,8,12,18\})}\big)
&= \big(1, (0, -1, -2), (-1, 0), (1, 2, 3), (1, -1, -3, -5), (-2, 1, 4), (3, 2), \\
&\qquad (-1, -4, -7), (-3, 1, 5, 9), 
(4, -1, -6), (-5, -4), (1, 6, 11),\dots \big).
\end{align*}
}
\end{exa}

By Proposition \ref{pro:reduction_rules}, we have the following result.

\begin{cor} \label{cor:V_odd_m_even}  For any even integer $m\ge 2$ and any set $V\subset [m]$ consisting of odd elements, we have $H_{n}\big(D^{(m,V)} \big)=1$ for every $n$.
\end{cor}

\section{Proof of Theorem \ref{thm:Main II}}

Given a primitive sequence $(t_1,\dots,t_{\ell-1})$ for some $\ell\ge 2$, let $(T_1,\dots,T_{\ell-1})$ be the dual sequence associated with $(t_1,\dots,t_{\ell-1})$, where $T_j=h(t_1,\dots,t_j)$. Recall that $T_{\ell-1}=0$ and $T_j\neq 0$ for all $j$, $1\le j\le \ell-2$. Let $V=\{2s,2(s+t_1),\dots,2(s+t_1+\cdots+t_{\ell-1})\}$ for some $s\ge 1$ and $m\ge\max(V)$. We shall prove that the sequence $H(D^{(m,V)})$ is periodic. Consider the two situations according to the parity of the modulus $m$.

%
%
\smallskip
Case 1. $m$ is even, say $m=2p$. By Theorem \ref{thm:Main I} along with $t_{\ell}=p-(t_1+\cdots+t_{\ell-1})\ge s$, the sequence $H\big(D^{(m,V)}\big)$ can be partitioned into sections of the form
\begin{equation} \label{eqn:even_and_odd_form}
\begin{aligned}
H\big(D^{(m,V)}\big) &=(\underbrace{1,\dots,1}_{s},\underbrace{a_{1,1},\dots,a_{1,t_1}}_{t_1},\underbrace{a_{2,1},\dots,a_{2,t_2}}_{t_2},\dots,\underbrace{a_{\ell,1},\dots,a_{\ell,t_{\ell}}}_{t_{\ell}},  \\
             & \qquad \underbrace{a_{\ell+1,1},\dots,a_{\ell+1,t_1}}_{t_1},\underbrace{a_{\ell+2,1},\dots,a_{\ell+2,t_2}}_{t_2},\dots,\underbrace{a_{2\ell,1},\dots,a_{2\ell,t_{\ell}}}_{t_{\ell}},\dots )
\end{aligned}
\end{equation}
such that each section $(a_{j,1}, \dots,a_{j,t_{\overline{j}}})$, $\overline{j}\equiv j\pmod{\ell}$, is an arithmetic progression with the common difference $d_j$, where
$a_{1,1}=0$, $d_1=-1$, $a_{j+1,1}= d_{j}$, and $d_{j+1} = d_{j}-a_{j,t_{\overline{j}}}$ for $j\ge 1$.  For $1\le j\le \ell-1$, we observe that 
\begin{equation} \label{eqn:recurrence}
d_{j+1}=d_j(2-t_j)-d_{j-1},
\end{equation}
where $d_0=a_{1,1}=0$. Notice that 
\[
\frac{d_{j+1}}{d_j}=2-t_j-\frac{d_{j-1}}{d_j},
\]
and that
\begin{equation} \label{eqn:product}
d_{j+1}=d_jT_j=d_1T_1\cdots T_j,
\end{equation}
for $1\le j\le \ell-1$.
It follows that $d_{\ell}=0$ and $d_j\neq 0$ for all $j$, $1\le j\le\ell-1$.

To prove the periodicity of $H\big(D^{(m,V)}\big)$, it suffices to prove that either $(d_{\ell-1},d_{\ell})=(1,0)$ or $(d_{2\ell-1},d_{2\ell})=(1,0)$ since in the former case, we have $a_{\ell,1}=\cdots=a_{\ell,t_{\ell}}=1$,  $a_{\ell+1,1}=0$, $d_{\ell+1}=-1$ and periodicity of $p$, while in the latter case, $a_{2\ell,1}=\cdots=a_{2\ell,t_{\ell}}=1$, $a_{2\ell+1,1}=0$, $d_{2\ell+1}=-1$ and periodicity of $2p$.

By Eq.(\ref{eqn:recurrence}), $\gcd(|d_{\ell-1}|,|d_{\ell-2}|)=\cdots=\gcd(|d_2|,|d_1|)=1$. Moreover, $d_{\ell-1}$ divides $d_{\ell-2}$ since $d_{\ell}=0$. Hence $d_{\ell-1}=1$ or  $-1$.

If $d_{\ell-1}=-1$ then by Eq.\,(\ref{eqn:product}), we have $T_1\cdots T_{\ell-2}=d_{\ell-1}/d_1=1$. Similar to Eq.\,(\ref{eqn:recurrence}), we have $d_{\ell+1}=d_{\ell}(2-t_{\ell})-d_{\ell-1}$ and
\[
d_{\ell+j+1}=d_{\ell+j}(2-t_j)-d_{\ell+j-1},
\]
for $1\le j\le \ell-1$. It follows that $d_{\ell+1}=-d_{\ell-1}=1$ and $d_{\ell+j+1}=d_{\ell+1}T_1\cdots T_j$ for $1\le j\le\ell-1$. Hence $d_{2\ell-1}=d_{\ell+1}T_1\cdots T_{\ell-2}=1$ and $d_{2\ell}=d_{2\ell-1}T_{\ell-1}=0$. This proves Case 1.

\smallskip
Case 2. $m$ is odd, say $m=2q-1$. Let $s=1$, i.e., $V=\{2,2(1+t_1),\dots,2(1+t_1+\cdots+t_{\ell-1})\}$, and let $t_{\ell}=m-(t_1+\cdots+t_{\ell-1})$. Notice that $2\in V-2(t_1+\cdots+t_i)$ for $1\le i\le \ell-1$. By the same argument as in Eqs.\,(\ref{eqn:initial_item_section}) and (\ref{eqn:first_run}), we obtain 
\begin{align*}
a_{b+1,k} &=H_{t_1+\cdots+t_b+k+1}\big(D^{(m,V)} \big) \\
          &=a_{b,t_b}H_{k+1}\big(D^{(m,V-2(t_1+\cdots+t_b))} \big)-d_bH_{k+1}\big(-1+D^{(m,V-2(t_1+\cdots+t_b))}\big) \\
\end{align*}
for $1\le b\le \ell-1$. Notice that in the evaluation of the section $(a_{\ell,1},\dots,a_{\ell,t_{\ell}})$, there is a concern that $1\in V-2q$ and $1\in V-2(q+t_1+\cdots+t_i)$ for $1\le i\le \ell-1$ and then the reduction rule in Proposition \ref{pro:reduction_rules}(ii) does not apply to the series $-1+D^{(m,V-2q)}$ and $-1+D^{(m,V-2(q+t_1+\cdots+t_i))}$. However, 
the entries of the sequences $H\big(-1+D^{(m,V-2q)}\big)$ and $H\big(-1+D^{(m,V-2(q+t_1+\cdots+t_i))}\big)$ do not appear in the evaluation since $a_{\ell-1,t_{\ell-1}}-d_{\ell-1}=d_{\ell}=0$, resulting from the primitive sequence $(t_1,\dots,t_{\ell-1})$, as shown below.
\begin{align*}
a_{\ell,k}&=H_{t_1+\cdots+t_{\ell-1}+k+1}\big(D^{(m,V)} \big) \\
      &=a_{\ell-1,t_{\ell-1}}H_{k+1}\big(D^{(m,V-2(t_1+\cdots+t_{\ell-1}))}\big)-d_{\ell-1}H_{k+1}\big(-1+D^{(m,V-2(t_1+\cdots+t_{\ell-1}))}\big) \\
      &=d_{\ell-1}H_{k}\big(D^{(m,V-2(t_1+\cdots+t_{\ell-1}+1))}\big)+(a_{\ell-1,t_{\ell-1}}-d_{\ell-1})H_{k}\big(-1+D^{(m,V-2(t_1+\cdots+t_{\ell-1}+1))}\big) \\
      &=d_{\ell-1}H_{k}\big(D^{(m,V-2(t_1+\cdots+t_{\ell-1}+1))}\big)  \\
      &= \cdots \\
      &=d_{\ell-1}H_{1}\big(D^{(m,V-2(t_1+\cdots+t_{\ell-1}+k))}\big) \\
      &=d_{\ell-1}
\end{align*}      
for $1\le k\le t_{\ell}$. That the sequence $H\big(D^{(m,V)}\big)$ is periodic with a period of either $m$ or $2m$ can be proved as in the previous case.
\qed

\begin{exa} {\rm
(i) Take $m=22$ and $V=\{2,8,12,14\}$, a feasible set derived from the primitive sequence $(3,2,1)$. Note that the periodicity of the sequence $H\big(D^{(22,\{2,8,12,14\})}\big)$ is 11,
\begin{align*}
H\big(D^{(22,\{2,8,12,14\})}\big) &=
\big(1, ((0, -1, -2), (-1, 0), (1), (1, 1, 1, 1, 1))^* \big) \\
&= (1, 0, -1, -2, -1, 0, 1, 1, 1, 1, 1)^*.
\end{align*}
(ii) If $m=21$ and $V=\{2,8,12,14\}$, the periodicity of the sequence $H\big(D^{(21,\{2,8,12,14\})}\big)$ is 21,
\begin{align*}
H\big(D^{(21,\{2,8,12,14\})}\big) &=
\big(1, ((0, -1, -2), (-1, 0), (1), (1, 1, 1, 1, 1, 1, 1, 1, 1, 1, 1, 1, 1, 1, 1))^* \big)\\
&= (1, 0, -1, -2, -1, 0, 1, 1, 1, 1, 1, 1, 1, 1, 1, 1, 1, 1, 1, 1, 1)^*
\end{align*}

(iii) Take $m=24$ and $V=\{2,8,16,22\}$, a feasible set derived from the primitive sequence $(3,4,3)$. Note that the periodicity of the sequence $H\big(D^{(24,\{2,8,16,22\})}\big)$ is 24,
\begin{align*}
H\big(D^{(24,\{2,8,16,22\})}\big) &=
\big(1, ((0, -1, -2), (-1, 0, 1, 2), (1, 0, -1), (-1, -1), \\
&\qquad (0, 1, 2), (1, 0, -1, -2), (-1, 0, 1), (1, 1))^*\big) \\
&=(1, 0, -1, -2, -1, 0, 1, 2, 1, 0, -1, -1, -1, 0, 1, 2, 1, 0, -1, -2, -1, 0, 1, 1)^*.
\end{align*}

(iv) If $m=23$ and $V=\{2,8,16,22\}$,  the periodicity of the sequence $H\big(D^{(23,\{2,8,16,22\})}\big)$ is 46,
\begin{align*}
H\big(D^{(23,\{2,8,16,22\})}\big) &=
\big(1, ((0, -1, -2), (-1, 0, 1, 2), (1, 0, -1), \\
&\qquad (-1, -1, -1, -1, -1, -1, -1, -1, -1, -1, -1, -1, -1), \\
&\qquad (0, 1, 2), (1, 0, -1, -2), (-1, 0, 1), (1, 1, 1, 1, 1, 1, 1, 1, 1, 1, 1, 1, 1))^*\big) \\
&=(1, 0, -1, -2, -1, 0, 1, 2, 1, 0,  \\
&\qquad -1, -1, -1, -1, -1, -1, -1, -1, -1, -1, -1, -1, -1, \\
&\qquad 0, 1, 2, 1, 0, -1, -2, -1, 0, 1, 1, 1, 1, 1, 1, 1, 1, 1, 1, 1, 1, 1)^*.
\end{align*}
}
\end{exa}

\section{Synthesis of feasible sets for periodic sequences}

In the section, we present a synthesis  of two feasible sets derived from primitive sequences (Propositions \ref{pro:synthesis_2s} and \ref{pro:synthesis_2s-1}), which can be generalized to Theorem \ref{thm:variety}.

\begin{pro} \label{pro:synthesis_2s}
For primitive sequences $\alpha=(t_1,\dots,t_{\ell_1-1})$ and $\beta=(r_1,\dots,r_{\ell_2-1})$, let $V_{\alpha}=\{2,2(1+t_1),\dots,2(1+t_1\cdots+t_{\ell_1-1})\}$ and $V_{\beta}=\{2,2(1+r_1),\dots,2(1+r_1+\cdots+r_{\ell_2-1})\}$. If $V=(V_{\alpha}+2s)\cup(V_{\beta}+k)$ for some $s\ge 0$,  $k\ge\max(V_{\alpha}+2s)-1$, and $m\ge\max(V)$ then the sequence of Hankel determinants of the series $D^{(m,V)}$ is periodic with a period of either $p$ or $2p$, where
\[
p=\begin{cases}
\frac{m}{2} &\mbox{for $m$ even} \\
m &\mbox{for $m$ odd.}
\end{cases}
\]
\end{pro}

\begin{proof} We consider the synthesis of $V$ according to the parity of $k$.

Case 1. $k$ is even. By Theorem \ref{thm:Main II}, the initial $p+s+1$ entries of the sequence $H\big(D^{(m,V)}\big)$ can be partitioned into sections of the form
\begin{equation} \label{eqn:synthesis_feasible_set}
\begin{aligned}
H\big(D^{(m,V)}\big) &=(\underbrace{1,\dots,1}_{s+1},\underbrace{a_{1,1},\dots,a_{1,t_1}}_{t_1},\underbrace{a_{2,1},\dots,a_{2,t_2}}_{t_2},\dots,\underbrace{a_{\ell_1,1},\dots,a_{\ell_1,t_{\ell_1}}}_{t_{\ell_1}},  \\
             & \qquad \underbrace{a_{\ell_1+1,1},\dots,a_{\ell_1+1,r_1}}_{r_1},\underbrace{a_{\ell_1+2,1},\dots,a_{\ell_1+2,r_2}}_{r_2},\dots,\underbrace{a_{\ell_1+\ell_2,1},\dots,a_{\ell_1+\ell_2,r_{\ell_2}}}_{r_{\ell_2}},\dots ),
\end{aligned}
\end{equation}
where $t_{\ell_1}=\frac{k+2}{2}-(s+1+t_1+\cdots+t_{\ell_1-1})$ and $r_{\ell_2}=p-\frac{k+2}{2}-(r_1+\cdots+r_{\ell_2-1})+s+1$.
Let $(T_1,\dots,T_{\ell_1-1})$ ($(R_1,\dots,R_{\ell_2-1})$, respectively) be the dual sequence associated with $\alpha$ ($\beta$, respectively), where $T_i=h(t_1,\dots,t_i)$ and $R_j=h(r_1,\dots,r_j)$. Recall that $T_{\ell_1-1}=0$ and $R_{\ell_2-1}=0$.
By the same argument as in Eqs.\,(\ref{eqn:recurrence}) and (\ref{eqn:product}), we have $d_{i+1}=d_1T_1\cdots T_i$ for $1\le i\le \ell_1-1$, and hence $d_{\ell_1}=0$ and $d_{\ell_1-1}=1$ or $-1$.

Let $q=\frac{k}{2}$. Notice that $2\in V-2q$ and $2\in V-2(q+r_1+\cdots+r_j)$ for $1\le j\le \ell_2-1$. Then
$a_{\ell_1+1,1}=d_{\ell_1}=0$, $d_{\ell_1+1}=-d_{\ell_1-1}$ and for $1\le j\le \ell_2-1$, we have
\begin{equation} \label{eqn:recurrence_2}
d_{\ell_1+j+1}=d_{\ell_1+j}(2-r_j)-d_{\ell_1+j-1},
\end{equation}
and hence $d_{\ell_1+j+1}=d_{\ell_1+1}R_1\cdots R_j$. It follows that $d_{\ell_1+\ell_2}=0$ and $d_{\ell_1+\ell_2-1}$ divides $d_{\ell_1+\ell_2-2}$. By Eq.,(\ref{eqn:recurrence_2}), $\gcd(|d_{\ell_1+\ell_2-1}|,|d_{\ell_1+\ell_2-2}|)=\cdots=\gcd(|d_{\ell_1+2}|,|d_{\ell_1+1}|)=1$. Hence $d_{\ell_1+\ell_2-1}=1$ or $-1$.  By the same argument as in the proof of Theorem \ref{thm:Main II}, the sequence $H\big(D^{(m,V)}\big)$ is periodic with a period of $p$ or $2p$.

Case 2. $k$ is odd. Notice that the set $V_{\beta}+k$ consists of odd elements.  If $m$ is even then $2\not\in (V_{\beta}+k)-2j$, modulo $m$, for all $j$.  Note that $2\in V-m-2s$ and $2\in V-m-2(s+t_1+\cdots+t_i)$ for $1\le i\le\ell_1-1$. We observe that the sequence  $H\big(D^{(m,V)}\big)$ can be partitioned into sections of the form
\begin{equation} \label{eqn:k_odd_m_even}
\begin{aligned}
H(D^{(m,V)}) &=(\underbrace{1,\dots,1}_{s+1},\underbrace{a_{1,1},\dots,a_{1,t_1}}_{t_1},\underbrace{a_{2,1},\dots,a_{2,t_2}}_{t_2},\dots,\underbrace{a_{\ell_1,1},\dots,a_{\ell_1,t_{\ell_1}}}_{t_{\ell_1}},  \\
             & \qquad \underbrace{a_{\ell_1+1,1},\dots,a_{\ell_1+1,t_1}}_{t_1},\underbrace{a_{\ell_1+2,1},\dots,a_{\ell_1+2,t_2}}_{t_2},\dots,\underbrace{a_{2\ell_1,1},\dots,a_{2\ell_1,t_{\ell_1}}}_{t_{\ell_1}},\dots ),
\end{aligned}
\end{equation}
where $t_{\ell_1}=\frac{m}{2}-(t_1+\cdots+t_{\ell_1-1})\ge s+1$. 
Moreover, if $m$ is odd then $2\in V-(m+k)$ and $2\in V-(m+k)-2(r_1+\cdots+r_j)$ for $1\le j\le \ell_2-1$. We observe that the initial $m+s+1$ entries of the sequence  $H\big(D^{(m,V)}\big)$ can be partitioned into sections of the form in Eq.\,(\ref{eqn:synthesis_feasible_set}), where $t_{\ell_1}=\frac{m+k+2}{2}-(s+1+t_1+\cdots+t_{\ell_1-1})$ and $r_{\ell_2}=m-\frac{m+k+2}{2}-(r_1+\cdots+r_{\ell_2-1})+s+1$.

The periodicity of the sequence $H\big(D^{(m,V)}\big)$ can be proved by the argument as in the previous case.
\end{proof}

\begin{exa} {\rm
(i) Take $m=22$, $V_{\alpha}=\{2,8,12,14\}$ and $V_{\beta}=\{2,6\}$, and  let $V=V_{\alpha}\cup(V_{\beta}+16)$. We observe that the periodicity of the sequence $H\big(D^{(22,\{2,8,12,14,18,22\})}\big)$ is 22, 
\begin{align*}
H\big(D^{(22,\{2,8,12,14,18,22\})}\big) 
&=\big(1, ((0, -1, -2), (-1, 0), (1), (1, 1), (0, -1), (-1), \\
&\qquad (0, 1, 2), (1, 0), (-1), (-1, -1), (0, 1), (1))^*\big) \\
&= (1, 0, -1, -2, -1, 0, 1, 1, 1, 0, -1, -1, 0, 1, 2, 1, 0, -1, -1, -1, 0, 1)^*.
\end{align*}

(ii) If $m=22$ and $V=V_{\alpha}\cup(V_{\beta}+15)$, the periodicity of $H\big(D^{(22,\{2,8,12,14,17,21\})}\big)$ is 11,
\begin{align*}
H\big(D^{(22,\{2,8,12,14,17,21\})}\big) 
&=
\big(1, ((0, -1, -2), (-1, 0), (1), (1, 1, 1, 1, 1))^*\big) \\
&=(1, 0, -1, -2, -1, 0, 1, 1, 1, 1, 1)^*.
\end{align*}

(iii) If $m=21$ and $V=V_{\alpha}\cup(V_{\beta}+15)$, the periodicity of $H\big(D^{(21,\{2,8,12,14,17,21\})}\big)$ is 42,
\begin{align*}
& H\big(D^{(21,\{2,8,12,14,17,21\})}\big) \\
&\quad= \big(1, ((0, -1, -2), (-1, 0), (1), (1, 1, 1, 1, 1, 1, 1, 1, 1, 1, 1, 1), (0, -1), (-1), \\
&\quad\qquad (0, 1, 2), (1, 0), (-1), (-1, -1, -1, -1, -1, -1, -1, -1, -1, -1, -1, -1), (0, 1), (1))^* \big) \\
&\quad= (1, 0, -1, -2, -1, 0, 1, 1, 1, 1, 1, 1, 1, 1, 1, 1, 1, 1, 1, 0, -1, -1, \\
&\quad\qquad 0, 1, 2, 1, 0, -1, -1, -1, -1, -1, -1, -1, -1, -1, -1, -1, -1, -1, 0, 1)^*.
\end{align*}
}
\end{exa}

\begin{pro} \label{pro:synthesis_2s-1}
For primitive sequences $\alpha=(t_1,\dots,t_{\ell_1-1})$ and $\beta=(r_1,\dots,r_{\ell_2-1})$, let $V_{\alpha}=\{2,2(1+t_1),\dots,2(1+t_1\cdots+t_{\ell_1-1})\}$ and $V_{\beta}=\{2,2(1+r_1),\dots,2(1+r_1+\cdots+r_{\ell_2-1})\}$. If $V=(V_{\alpha}+2s-1)\cup(V_{\beta}+k)$ for some $s\ge 0$ and  $k\ge\max(V_{\alpha}+2s-1)-1$ then the sequence of Hankel determinants of the series $D^{(m,V)}$ is periodic for every $m\ge\max(V)$.
\end{pro}

\begin{proof} Following the proof of Proposition \ref{pro:synthesis_2s}, it suffices to consider the reductions when the element 2 appears in the avoiding set. Notice that $V_{\alpha}+2s-1$ consists of odd elements. Consider the synthesis of $V$ according to the parity of $k$.

Case 1. $k$ is even. If $m$ is even then $2\not\in (V_{\alpha}+2s-1)-2j$ for all $j$. Notice that $2\in V-k$ and $2\in V-k-2(r_1+\cdots+r_j)$ for $1\le j\le \ell_2-1$. By the same argument as in the proof of Proposition \ref{pro:synthesis_2s}, we observe that the sequence $H\big(D^{(m,V)}\big)$ can be partitioned into sections of the form
\begin{equation}
\begin{aligned}
H(D^{(m,V)}) &=(\underbrace{1,\dots,1}_{(k+2)/2},\underbrace{a_{1,1},\dots,a_{1,r_1}}_{r_1},\underbrace{a_{2,1},\dots,a_{2,r_2}}_{r_2},\dots,\underbrace{a_{\ell_2,1},\dots,a_{\ell_2,r_{\ell_2}}}_{r_{\ell_2}},  \\
             & \qquad \underbrace{a_{\ell_2+1,1},\dots,a_{\ell_2+1,r_1}}_{r_1},\underbrace{a_{\ell_2+2,1},\dots,a_{\ell_2+2,r_2}}_{r_2},\dots,\underbrace{a_{2\ell_2,1},\dots,a_{2\ell_2,r_{\ell_2}}}_{r_{\ell_2}},\dots ),
\end{aligned}
\end{equation}
where $r_{\ell_2}=\frac{m}{2}-(r_1+\cdots+r_{\ell_2-1})$. Moreover, if $m$ is odd then $2\in V-(m-2s-1)$ and $2\in V-(m-2s-1)-2(t_1+\cdots+t_i)$ for $1\le i\le \ell_1-1$. The initial entries of the sequence $H\big(D^{(m,V)}\big)$ can be partitioned into sections of the form
\begin{equation}
\begin{aligned}
H(D^{(m,V)}) &=(\underbrace{1,\dots,1}_{(k+2)/2},\underbrace{a_{1,1},\dots,a_{1,r_1}}_{r_1},\underbrace{a_{2,1},\dots,a_{2,r_2}}_{r_2},\dots,\underbrace{a_{\ell_2,1},\dots,a_{\ell_2,r_{\ell_2}}}_{r_{\ell_2}},  \\
             & \qquad \underbrace{a_{\ell_2+1,1},\dots,a_{\ell_2+1,t_1}}_{t_1},\underbrace{a_{\ell_2+2,1},\dots,a_{\ell_2+2,t_2}}_{t_2},\dots,\underbrace{a_{\ell_2+\ell_1,1},\dots,a_{\ell_2+\ell_1,t_{\ell_1}}}_{t_{\ell_1}},\dots ),
\end{aligned}
\end{equation}
where $r_{\ell_2}=\frac{m-k-1}{2}-(r_1+\cdots+r_{\ell_2-1})+s$ and $t_{\ell_1}=\frac{m+k+1}{2}-(t_1+\cdots+t_{\ell_1-1})-s$.

Case 2. $k$ is odd. Notice that $V$ consists entirely of odd elements. If $m$ is even then by Corollary \ref{cor:V_odd_m_even}, we have $H_n(D^{(m,V)})=1$ for all $n$. The periodicity is 1. Moreover, if $m$ is odd then $2\in V-(m+2s-1)$ and $2\in V-(m+2s-1)-2(t_1+\cdots+t_i)$ for $1\le i\le \ell_1-1$. Moreover, $2\in V-(m+k)$ and $2\in V-(m+k)-(r_1+\cdots+r_j)$ for for $1\le j\le \ell_2-1$. The initial entries of the sequence $H\big(D^{(m,V)}\big)$ can be partitioned into sections of the form
\begin{equation} 
\begin{aligned}
H\big(D^{(m,V)}\big) &=(\underbrace{1,1,\dots,1}_{(m+2s+1)/2},\underbrace{a_{1,1},\dots,a_{1,t_1}}_{t_1},\underbrace{a_{2,1},\dots,a_{2,t_2}}_{t_2},\dots,\underbrace{a_{\ell_1,1},\dots,a_{\ell_1,t_{\ell_1}}}_{t_{\ell_1}},  \\
             & \qquad \underbrace{a_{\ell_1+1,1},\dots,a_{\ell_1+1,r_1}}_{r_1},\underbrace{a_{\ell_1+2,1},\dots,a_{\ell_1+2,r_2}}_{r_2},\dots,\underbrace{a_{\ell_1+\ell_2,1},\dots,a_{\ell_1+\ell_2,r_{\ell_2}}}_{r_{\ell_2}},\dots ),
\end{aligned}
\end{equation}
where $t_{\ell_1}=\frac{k+1}{2}-(s+t_1+\cdots+t_{\ell_1-1})$ and $r_{\ell_2}=m-\frac{k+1}{2}-(r_1+\cdots+r_{\ell_2-1})+s$.
The periodicities of the sequences $H\big(D^{(m,V)}\big)$ can be obtained as in the proof of Proposition \ref{pro:synthesis_2s}.
\end{proof}

\smallskip
Making use of the same argument as in the proofs of Propositions \ref{pro:synthesis_2s} and \ref{pro:synthesis_2s-1}, we can prove Theorem \ref{thm:variety} by induction.

\begin{exa} {\rm
Take $m=21$, $V_{\alpha}=\{2,8,12,14\}$ and $V_{\beta}=\{2,6\}$, and  let $V=(V_{\alpha}-1)\cup(V_{\beta}+14)$. We observe that the periodicity of the sequence $H\big(D^{(21,\{1,7,11,13,16,20\})}\big)$ is 42, 
\begin{align*}
& H\big(D^{(21,\{1,7,11,13,16,20\})}\big) \\
&\quad =\big(1, 1, 1, 1, 1, 1, 1, 1, ((0, -1), (-1), (0, 1, 2), (1, 0), (-1), (-1, -1, -1, -1, -1, -1, -1, -1\\ 
&\quad\qquad  -1, -1, -1, -1), (0, 1), (1), (0, -1, -2), (-1, 0), (1), (1, 1, 1, 1, 1, 1, 1, 1, 1, 1, 1, 1))^* \big) \\
&\quad = (1, 1, 1, 1, 1, 1, 1, 1, 0, -1, -1, 0, 1, 2, 1, 0, -1, \\
&\quad\qquad -1, -1, -1, -1, -1, -1, -1, -1, -1, -1, -1, -1, 0, 1, 1, 0, -1, -2, -1, 0, 1, 1, 1, 1, 1)^*.
\end{align*}
}
\end{exa}

In the following, we show an extension of admissible sequences and a construction for primitive sequences from admissible sequences. 

\begin{thm} Given an admissible sequence $(t_1,\dots,t_b)$, $b\ge 1$, let $(T_1,\dots,T_b)$ be the dual sequence associated to $(t_1,\dots,t_b)$, where $T_j=h(t_1,\dots,t_j)$. Then
the following properties hold.
\begin{enumerate}
\item  For $1 \le j\le b$, $\prod_{k=1}^j T_k$ is an integer.

\item If $|\prod_{k=1}^b T_k|=1$ and either $T_b=1$ or $T_b<0$ then there exists a positive integer $t_{b+1}$ such that the sequence $(t_1,\dots,t_{b+1})$ is primitive.

\item If $|\prod_{k=1}^b T_k|=1$ and either $T_b=\frac{1}{2}$ or $T_b<0$ then there exists a positive integer $t_{b+1}$ such that  $T_{b+1}=-1$.

\item If $|\prod_{k=1}^b T_k|=1$ and $T_b<0$ then there exists a positive integer $t_{b+1}$ such that  $T_{b+1}=1$.

\end{enumerate} 
\end{thm}

\begin{proof} (i) By the definition of admissible  sequence, notice that
\begin{equation} \label{eqn:T_j+1}
T_{j+1}=2-t_{j+1}-\frac{1}{T_j},
\end{equation}
and that both $T_1=2-t_1$ and $T_1T_2=T_1(2-t_2)-1$ are integers. Multiplying both sides of Eq.\,(\ref{eqn:T_j+1}) by $T_1\cdots T_j$ yields
\[
\prod_{k=1}^{j+1} T_k=(2-t_{j+1})\prod_{k=1}^{j} T_k-\prod_{k=1}^{j-1}T_k,
\] which is an integer by induction. The assertion follows.

\smallskip
(ii) If $|\prod_{k=1}^b T_k|=1$ then by (i), $|T_b|^{-1}=|\prod_{k=1}^{b-1} T_k|$ is an integer. By Eq.\,(\ref{eqn:T_j+1}) along with $T_{b+1}=0$, we have $r_{b+1}=2-T_b^{-1}\ge 1$.

\smallskip
(iii) Notice that $T_b^{-1}$ is either 2 or a negative integer. By Eq.\,(\ref{eqn:T_j+1}) along with $T_{b+1}=-1$, we have $r_{b+1}=2-T_{b+1}-T_b^{-1}\ge 1$.

\smallskip
(iv) Notice that $T_b^{-1}$ is a negative integer. By Eq.\,(\ref{eqn:T_j+1}) along with $T_{b+1}=1$, we have $r_{b+1}=2-T_{b+1}-T_b^{-1}\ge 2$.
\end{proof}

\begin{exa} {\rm (i) Given the admissible sequence $(t_1,t_2,t_3)=(5,3,4)$, we have $(T_1,T_2,T_3)=(-3,-\frac{2}{3},-\frac{1}{2})$.  To construct a primitive sequence $(t_1,t_2,t_3,t_4)$, i.e., for $T_4=0$, we take $t_4=4$. 

(ii) If $(t_1,t_2,t_3)=(5,3,3)$ then $(T_1,T_2,T_3)=(-3,-\frac{2}{3},\frac{1}{2})$.  To construct an admissible sequence $(t_1,t_2,t_3,t_4)$ with $T_4=-1$, we take $t_4=1$. 
}
\end{exa}

\section{Some sporadic cases of periodic sequences}

There are still series $D^{(m,V)}$ with periodic sequences $H\big(D^{(m,V)} \big)$ not covered by Theorem \ref{thm:Main II} or Theorem \ref{thm:variety}. In the following, we present some cases bypassing the lack of reduction rule for $H_n\big(-1+D^{(m,V)}\big)$ when $V$ contains the element 1.

\begin{pro} \label{pro:modulo_3} The sequences of Hankel determinants of the series $D^{(3,\{1\})}$, $D^{(3,\{2\})}$ and $D^{(3,\{3\})}$ are periodic with a period of 10, 
\begin{align*}
H\big(D^{(3,\{1\})} \big) &=(1,1,0,-1,-1,-1,-1,0,1,1)^*, \\
H\big(D^{(3,\{2\})} \big) &=(1,0,-1,-1,-1,-1,0,1,1,1)^*, \\
H\big(D^{(3,\{3\})} \big) &=(1,1,1,0,-1,-1,-1,-1,0,1)^*. \\
\end{align*}
\end{pro}

\begin{proof}  By Proposition \ref{pro:reduction_rules}(i), we have 
\[
H_{n+2}\big(D^{(3,\{3\})}\big)=H_{n+1}\big(D^{(3,\{1\})}\big)=H_n\big(D^{(3,\{2\})}\big).
\] 
The relation $D^{(3,\{2\})}=x+D^{(3,\{1\})}$ can be obtained by solving the equations in Eq.\,(\ref{eqn:continued_fraction}) for the series $D^{(3,\{2\})}$ and $D^{(3,\{1\})}$.
By Proposition \ref{pro:reduction_rules}(i) and Lemma \ref{lem:-H(2)}, we have
\begin{align*}
H_n\big(D^{(3,\{2\})} \big) &= H_{n-1}\big(-1+D^{(3,\{3\})} \big) \\
               &= -H_{n-3}^{(2)}\big(D^{(3,\{1\})} \big) \\
               &= -H_{n-3}^{(2)}\big(-x+D^{(3,\{2\})} \big).
\end{align*}
Notice that $H_{n-3}^{(2)}\big(-x+D^{(3,\{2\})} \big)=H_{n-3}^{(2)}\big(D^{(3,\{2\})} \big)$ since these two series have the same coefficient of $x^j$ for all $j\ge 2$. Hence by Lemma \ref{lem:H(2)}, we have
\begin{align*}
H_n\big(D^{(3,\{2\})} \big) &= -H_{n-3}^{(2)}\big(D^{(3,\{2\})} \big) \\
               &= -H_{n-2}\big(D^{(3,\{2\})} \big)+H_{n-2}\big(-1+D^{(3,\{2\})} \big) \\
               &= -H_{n-3}\big(D^{(3,\{3\})} \big) \\
               &= -H_{n-4}\big(D^{(3,\{1\})} \big) \\
               &= -H_{n-5}\big(D^{(3,\{2\})} \big).
\end{align*}
Using Eq.\,(\ref{eqn:trivial}), the initial elements can be determined by the above reduction. It follows that
\[
H\big(D^{(3,\{2\})} \big)=(1,0,-1,-1,-1,-1,0,1,1,1)^*.
\]
\end{proof}

\begin{cor} We have $H_{n+1}\big(D^{(3,\{1,3\})} \big) =H_n\big(D^{(3,\{1,2\})} \big)$ and
\[
H\big(D^{(3,\{1,2\})} \big) =(1,0,-1,-1,-1,0,1,1)^*.
\]
\end{cor}

\begin{proof} Along with the relation $D^{(3,\{2,3\})}=x+D^{(3,\{1,3\})}$ derived from Eq.\,(\ref{eqn:continued_fraction}), the result $H_n\big(D^{(3,\{1,2\})}\big)=-H_{n-4}\big(D^{(3,\{1,2\})}\big)$ can be proved  by the same arguments as in the proof of Proposition \ref{pro:modulo_3}.
\end{proof}

\begin{pro} We have $H_{n+2}\big(D^{(5,\{1,5\})}\big)=H_{n+1}\big(D^{(5,\{3,4\})}\big)=H_n\big(D^{(5,\{1,2\})}\big)$,  and 
\[
H\big(D^{(5,\{1,2\})} \big)=(1,0,-1,-1,-1,-1,0,1,1,1)^*. 
\]
\end{pro}

\begin{proof} By Proposition \ref{pro:reduction_rules}, we have
\begin{align*}
H_n\big(D^{(5,\{1,2\})} \big) &= H_{n-1}\big(-1+D^{(5,\{4,5\})} \big) \\
               &= -H_{n-2}\big(D^{(5,\{2,3\})} \big)+H_{n-2}\big(-1+D^{(5,\{2,3\})} \big) \\
               &= -H_{n-3}\big(D^{(5,\{1,5\})} \big) \\
               &= -H_{n-4}\big(D^{(5,\{3,4\})} \big) \\
               &= -H_{n-5}\big(D^{(5,\{1,2\})} \big).
\end{align*}
The assertion follows.
\end{proof}

\section{Concluding remarks}
In this paper, we study the Hankel determinants for the  Dyck paths with peaks avoiding the heights in multiple congruence classes of a modulo $m$, which contain numerous instances with periodic sequences of Hankel determinants. 
One of our contributions is a sufficient condition for the avoiding set $(m,V)$ such that the sequence $H\big(D^{(m,V)}\big)$ is periodic, and a  construction for a variety of such sets.

For evaluating the Hankel determinants, we develop a unified approach by establishing fundamental reduction rules for the determinants. The pivotal points in the evaluation are the reductions when the  element 2 appears in the avoiding set.
An obstacle occurs in the process when $1\in V-2j$ for some $j$ since no deterministic rule is available for further reductions. Evidences shows that non-periodic sequences $H\big(D^{(m,V)}\big)$ tend to chaos out of control except for the sets $V$ consisting of even elements of an even modulo $m$. 
However, there are still a lot of instances with periodicity not covered by the sufficient condition. We raise a question about the periodicity of the sequence $H\big(D^{(m,V)}\big)$ for the Dyck paths with extremal avoiding sets $(m,V)$. An ultimate question is to characterize the set $(m,V)$ with periodic sequence $H\big(D^{(m,V)}\big)$.

\begin{con} For any integer $m\ge 3$ and the set $V=\{1,2,\dots,m-1\}$, the sequence of Hankel determinants of the series $D^{(m,V)}$ is periodic of the form
\[
H\big(D^{(m,V)}\big)=\begin{cases}
           (1,\underbrace{0,\dots, 0}_{m-2},1,1)^* & \mbox{if $m\equiv 1, 2 \pmod{4} $} \\
           (1,\underbrace{0,\dots, 0}_{m-2},-1,-1,-1,\underbrace{0,\dots, 0}_{m-2}, 1,1)^* & \mbox{if $m\equiv 0, 3 \pmod{4}$.}
  \end{cases}
\]

\end{con}

\section*{Acknowledgements.}
The authors thank the referees for reading the
manuscript carefully and providing helpful suggestions. 
The authors were supported in part by
Ministry of Science and Technology (MOST) grant 110-2115-M-003-011-MY3 (S.-P. Eu), and 109-2115-M-153-004-MY2 (T.-S. Fu).

\begin{table}[ht]
\caption{The periodicity of the sequences $H(D^{(m,V)})$ for all $(m,V)$ with $m\le 5$.}
{\small
\begin{tabular}{cp{4.2in}c}
\hline
$(m,V)$ & \multicolumn{1}{c}{$H(D^{(m,V)})$} &  period\\
\hline
$(2,\{1\})$ &  $(1)^*$ &  $1$\\
$(2,\{2\})$ &  $(1,0,-1,-1,0,1)^*$  &  $6$\\
\hline
$(3,\{1\})$ &  $(1,1,0,-1,-1,-1,-1,0,1,1)^*$ &  $10$\\
$(3,\{2\})$ &  $(1,0,-1,-1,-1,-1,0,1,1,1)^*$ &  $10$\\
$(3,\{3\})$ &  $(1,1,1,0,-1,-1,-1,-1,0,1)^*$ &  $10$\\
$(3,\{1,2\})$ &  $(1,0,-1,-1,-1,0,1,1)^*$ &  $8$\\
$(3,\{1,3\})$ &  $(1,1,0,-1,-1,-1,0,1)^*$ &  $8$\\
$(3,\{2,3\})$ &  $(1,0,0,-1,-1,0,0,1)^*$ &  $8$\\
\hline
$(4,\{1\})$ &  $(1)^*$ &  $1$\\
$(4,\{2\})$ &  $(1,0,-1,-1,-1,0,1,1)$ &  $8$\\
$(4,\{3\})$ &  $(1)^*$ &  $1$\\
$(4,\{4\})$ &  $(1,1,0,-1,-1,-1,0,1)^*$ &  $8$\\
$(4,\{1,2\})$ &  $(1,0,-1,0,1,1,1,0,0,-1,-1,-1,0,1,0,-1,-1,-1,0,0,1,1)^*$ &  $22$ \\
$(4,\{1,3\})$ &  $(1)^*$ &  $1$\\
$(4,\{1,4\})$ &  $(1,1,0,0,-1,-1,-1,0,1,0,-1,-1,-1,0,0,1,1,1,0,-1,0,1)^*$ &  $22$\\
$(4,\{2,3\})$ &  $(1,0,0,-1,-1,-1,0,1,0,-1,-1,-1,0,0,1,1,1,0,-1,0,1,1)^*$ &  $22$\\
$(4,\{2,4\})$ &  $(1,0,-1,-1,0,1)^*$ &  $6$\\
$(4,\{3,4\})$ & $(1,1,0,-1,0,1,1,1,0,0,-1,-1-1,0,1,0,-1,-1,-1,0,0,1)^*$ &  $22$\\
$(4,\{1,2,3\})$ &  $(1,0,0,-1,-1,-1,0,0,1,1)^*$ &  $10$\\
$(4,\{1,2,4\})$ &  $(1,0,-1,0,1)^*$ &  $5$\\
$(4,\{1,3,4\})$ &  $(1,1,0,0,-1,-1,-1,0,0,1)^*$ &  $10$\\
$(4,\{2,3,4\})$ &  $(1,0,0,0,1)^*$ &  $5$\\
\hline
$(5,\{1\})$ & $(1,1,1,0,-1,-2,-2,-3,-4,-5,-1,7,23,31,51,116,149,\cdots)$  & none\\
$(5,\{2\})$ & $(1,0,-1,-2,-2,-3,-4,-5,-1,7,23,31,51,116,149,118,-426,\cdots)$  & none\\
$(5,\{3\})$ & $(1,1,1,1,0,-1,-2,-2,-3,-4,-5,-1,7,23,31,51,116,149,\cdots)$  & none\\
$(5,\{4\})$ & $(1,1,0,-1,-2,-2,-3,-4,-5,-1,7,23,31,51,116,149,118,\cdots)$  & none\\
$(5,\{5\})$ & $(1,1,1,1,1,0,-1,-2,-2,-3,-4,-5,-1,7,23,31,51,116,149,\cdots)$  & none\\
$(5,\{1,2\})$ & $(1,0,-1,-1,-1,-1,0,1,1,1)^*$  & $10$\\
$(5,\{1,3\})$ & $(1,1,1,0,-1,-2,-2,-3,-4,-1,7,
15,23,47,68,53,-202,-618,\cdots)$  & none\\
$(5,\{1,4\})$ & $(1,1,0,-1,-2,-2,-3,-4,-1,7,15,23,47,68,53,-202,-618,\cdots$)  & none\\
$(5,\{1,5\})$ & $(1,1,1,0,-1,-1,-1,-1,0,1)^*$  & $10$\\
$(5,\{2,3\})$ & $(1,0,0,-1,-1,-1,0,0,1,1)^*$  & $10$\\
$(5,\{2,4\}$ & $(1,0,-1,-2,-2,-3,-4,-1,7,15,23,47,68,53,-202,-618,\cdots)$ & none\\
$(5,\{2,5\})$ & $(1,0,-1,-1,-2,-4,-2,5,13,20,43,67,60,-187,-595,-1338,\cdots)$ & none\\
$(5,\{3,4\})$ & $(1,1,0,-1,-1,-1,-1,0,1,1)^*$  & $10$\\
$(5,\{3,5\})$ & $(1,1,1,1,0,-1,-2,-2,-3,-4,-1,7,15,23,47,68,53,-202,\cdots)$ & none\\
$(5,\{4,5\})$ & $(1,1,0,0,-1,-1,-1,0,0,1)^*$ & $10$\\
$(5,\{1,2,3\})$ & $(1,0,0,-1,0,1,1,0,-1,0,0,1,1,1,0,0,0,1,1)^*$ & $19$\\
$(5,\{1,2,4\})$ & $(1,0,-1,-1,-1,-1,0,1,1,1)^*$  & $10$\\
$(5,\{1,2,5\})$ & $(1,0,-1,0,0,1,1,1,0,0,0,1,1,1,0,0,-1,0,1)^*$ & $19$\\
$(5,\{1,3,4\})$ & $(1,1,0,-1,-1,-1,-1,0,1,1)^*$ & $10$\\
$(5,\{1,3,5\})$ & $(1,1,1,0,-1,-1,-1,-1,0,1)^*$ & $10$\\
$(5,\{1,4,5\})$ & $(1,1,0,0,0,1,1,1,0,0,-1,0,1,1,0,-1,0,0,1)^*$ & $19$\\
$(5,\{2,3,4\})$ & $(1,0,0,0,1,1,1,0,0,-1,0,1,1,0,-1,0,0,1,1)^*$  & $19$\\
$(5,\{2,3,5\})$ & $(1,0,0,-1,-1,0,1,1,2,1,-1,-2,-2,-3,-1,2,3,3,4,1,-3,-4\cdots)$ & none\\
$(5,\{2,4,5\})$ & $(1,0,-1,-1,-2,-1,1,2,2,3,1,-2,-3,-3,-4,-1,3,4,4,5,1\cdots)$ & none\\
$(5,\{3,4,5\})$ & $(1,1,0,0,-1,0,1,1,0,-1,0,0,1,1,1,0,0,0,1)^*$ & $19$\\
$(5,\{1,2,3,4\})$ & $(1,0,0,0,1,1)^*$  & $6$\\
$(5,\{1,2,3,5\})$ & $(1,0,0,-1,0,1)^*$  & $6$\\
$(5,\{1,2,4,5\})$ & $(1,0,-1,0,0,1)^*$  & $6$\\
$(5,\{1,3,4,5\})$ & $(1,1,0,0,0,1)^*$  & $6$\\
$(5,\{2,3,4,5\})$ & $(1,0,0,0,0,1)^*$  & $6$\\
\hline
\end{tabular}
}
\label{tab:list}
\end{table}

\end{document}